\documentclass[a4paper,10pt]{amsart}
\usepackage{amssymb}
\usepackage{cmap}
\usepackage{tikz-cd}
\usepackage[pdfdisplaydoctitle=true,
            colorlinks=true,
            urlcolor=blue,
            citecolor=blue,
            linkcolor=blue,
            pdfstartview=FitH,
            pdfpagemode= UseNone,
            bookmarksnumbered=true,
            backref=page]{hyperref}
\usepackage{hyperxmp}
\usepackage{todonotes}
\usepackage[all]{xy}
\usepackage{ulem}
\usepackage{cancel}

\theoremstyle{plain}
\newtheorem{thm}{Theorem}[section]
\newtheorem*{thm*}{Theorem}

\newtheorem{lem}[thm]{Lemma}
\newtheorem{prop}[thm]{Proposition}

\theoremstyle{definition}

\theoremstyle{definition}
\newtheorem{rem}[thm]{Remark}
\newtheorem{defi}[thm]{Definition}

\newtheorem{expl}[thm]{Example}
\newtheorem*{ass*}{Assumption}

\numberwithin{equation}{section}

\allowdisplaybreaks

\usepackage{todonotes}
\hypersetup{
    pdftitle={Sobolev metrics on spaces of manifold valued curves}, %
    pdfsubject={MSC 2010: 58B20, 58D10, 35G55, 35G60}, %
    pdfkeywords={Sobolev metrics, manifold valued curves}, %
    pdflang=EN %
    }
\author{Yuxiu Lu }
\address{Department of Mathematics and statistics, Nanjing University of Science and Technology, 210094, P.R.China
 }
\email{luyxnj@gmail.com }
\date{September 15, 2026}



\begin{document}

\title[Vanishing geodesic distance phenomenon and boundedness of curvatures]{Vanishing geodesic distance phenomenon and boundedness of curvatures}

%
%

%
\begin{abstract}
This article concerns a conjecture proposed by Michor and Mumford in their 2005 paper, in which the vanishing of the geodesic distance is associated to the local unboundedness of the sectional curvature. In this article, we consider a Riemannian Banach manifold with all of its curvature tensors locally bounded. We derive several properties regarding the local boundedness of curvature tensors, and give a sufficient condition for a non-degenerate geodesic distance, which partially verifies the conjecture in its converse form.
\end{abstract}

\maketitle

\setcounter{tocdepth}{2}
\tableofcontents

\section{Introduction}
In the paper~\cite{michor2005}, Michor and Mumford found that the $L^2$-metric on the space of immersions  
$\operatorname{Imm}(M,N)$ induces vanishing geodesic distance on its quotient space
$\operatorname{Imm}(M,N)/\operatorname{Diff}(M)$. After some calculations, they also showed that "the curvature is positive and unbounded in some high
frequency directions", so they were in a position to conjecture a relationship between the vanishing of the geodesic distance and the local unboundedness of sectional curvatures. 

After this conjecture was proposed, considerable efforts have been made to understand the degeneracy of the geodesic distance on different diffeomorphism groups endowed with Sobolev metrics, as shown in~\cite[Section 1.2]{bauer2020}. But the reasons for this degeneracy remain mysterious. Recently, Magnani and Tiberio proved the conjecture for Hilbertian H-type groups equipped with weak left invariant Riemannian metrics; see~\cite{magnani2025}.
They showed the vanishing of geodesic distance that is related to the sub-Riemannian structure of the H-type groups, and found a sequence of $2$-dimensional subspaces such that their sectional curvatures go to infinity.

In the present article, we attack the problem from the other side. It is natural to consider an arbitrary infinite dimensional Riemannian manifold $M$ with all curvature tensors bounded. If the non-degeneracy of the geodesic distance is proved for such a manifold $M$, we essentially prove the equivalent form of the conjecture.

However, even the boundedness of curvature tensors appears to be more complicated than one could imagine, so it should be carefully discussed in Section~\ref{sec:boundcurv} and Appendix~\ref{sec:appendix}. It should also be noted that all of our discussions presuppose the existence of curvature tensors, and evidently presuppose  the existence of Levi-Civita connection and Christoffel symbols as well. Moreover, we restrict ourselves to a Banach manifold for its nice analogy to a finite-dimensional manifold.

The main theorem of this article asserts  that the local boundedness with an extra assumption on the existence of a local geodesic ball leads to the non-degeneracy of geodesic distance, which is proved step by step from Section~\ref{sec:geodesicball} to Section~\ref{sec:vanishing}. In fact, we show that within a geodesic ball, the uniform bound for curvatures gives a lower bound for the distance of any two conjugate points by an estimate on the length of Jacobi fields. This enables us to show that the exponential map must be a local diffeomorphism in the ball, which implies the non-degeneracy of geodesic distance by a well-known argument from basic Riemannian geometry in~\cite{docarmo1992}.

\section{Boundedness of curvature tensors }\label{sec:boundcurv}
In this section, we present several basic properties pertaining to boundedness of curvatures in infinite-dimensional setting. Note that in infinite-dimensional case, the Levi-Civita connection and the Christoffel symbols may not exist. But if we assume the existence of curvature tensors, it is automatically true that 
both the Levi-Civita connection and the Christoffel symbols exist.

Firstly we define the boundedness of curvatures with respect to a unique point or a local neighborhood.  

\begin{defi}[Boundedness of curvature]
A  Riemmannian manifold $(M,g)$ is said to have  bounded curvature at $x\in M$ if 
there exists a constant $N>0$ such that the curvature tensor $||R(X,Y)Z||_g\leq N$ for any unit vectors $X, Y, Z\in T_x M$. Otherwise $(M,g)$ is said to have a unbounded curvature at $x\in M$, i.e., if for any $N>0$ there exists 
unit vectors $X_N, Y_N, Z_N\in T_x M$ such that $||R(X_N,Y_N)Z_N||_g>N$. The  Riemmannian manifold $(M,g)$
is said to  have  locally bounded curvature if for any point $x\in M$ there exists a neighborhood $U_x$ with respect to the manifold topology  such that the curvature is uniformly bounded on $U_x$.
\end{defi}
In the following result, we  will see that the boundedness of the curvature tensors is equivalent to the the boundedness of all sectional curvatures:
\begin{prop}\label{seccurv}
For any normalized vectors $X, Y, Z$ at any point $p\in U_x$, the curvature tensors $R(X,Y)Z$ are bounded  if and only if the sectional curvatures $\langle R(X,Y)X, Y\rangle$ are bounded.   
\end{prop} 
\begin{proof}
First note that the curvature tensors $R(X,Y)Z$ are bounded for any normalized vectors $X, Y, Z$ if and only if $\langle R(X,Y)Z, W\rangle$ are bounded for any normalized vectors $X,Y,Z,W$. Therefore, it suffices to show that $\langle R(X,Y)Z, W\rangle$ are bounded given that $\langle R(X,Y)X, Y\rangle$ are bounded by a constant $C>0$.
By the linearity of curvature tensors, we have
\begin{align*}
\langle R(X+Z,Y)X+Z, Y\rangle=
\langle R(X,Y)X, Y\rangle+2\langle R(X,Y)Z, Y\rangle+\langle R(Z,Y)Z, Y\rangle.  \end{align*}
It follows by the triangle inequality and 
$||X+Z||\leq2$ that $||\langle R(X,Y)Z, Y\rangle||_g\leq3C$ for arbitrary choice of normalized vectors $X,Y,Z$.

Similarly we also have 
\begin{align*}
\langle R(X,Y+T)Z, Y+T\rangle&=\langle R(X,Y)Z, T\rangle+
\langle R(X,T)Z, Y\rangle\\&+\langle R(X,Y)Z, Y\rangle+\langle R(X,T)Z, T\rangle.
\end{align*}
Again by the triangle inequality, 
it follows that 
\begin{align*}
||\langle R(X,Y)Z, T\rangle-\langle R(Y,Z)X, T\rangle||_g=||\langle R(X,Y)Z, T\rangle+\langle R(X,T)Z, Y\rangle||_g
\leq9C.
\end{align*}
By a cyclic permutation of the first three elements, we obtain another inequality
\begin{align*}
||\langle R(Z,X)Y, T\rangle-\langle R(Y,Z)X, T\rangle||_g=||\langle R(Y,Z)X, T\rangle-\langle R(Z,X)Y, T\rangle||_g
\leq9C.
\end{align*}
Finally we add up the two inequalities and use the Bianchi identity
\begin{align*}
3||\langle R(Y,Z)X, T\rangle||_g&=||\langle R(Z,X)Y,T\rangle -\langle R(Y,Z)X, T\rangle\\&+\langle R(X,Y)Z, T\rangle-\langle R(Y,Z)X, T\rangle||_g
\leq18C.
\end{align*}
The assertion follows by the arbitrariness of $X,Y,Z,T$.
\end{proof}
Moreover for $M$ being a Banach manifold, the boundedness of curvature tensors is not pertaining to the whole curvature tensor.
\begin{prop}\label{wholecuraturetensor}
Let $M$ be a Banach manifold and let $X, Y, Z\in T_pM$ be unit vectors. If all  curvature tensors $R(X,Y)Z$ are bounded when two of the vectors $X, Y, Z$ are fixed, then all the curvature tensors are bounded.   
\end{prop}
\begin{proof}
Since the curvature tensor $R(X,Y)Z$ is linear with respect to any of $X, Y, Z$, without loss of generality, it suffices to show the proposition for $X$ varying. Therefore by assumption, we can define a collection of bounded linear operators $T$ via $T(X):=R(X,Y)Z$ for any fixed $Y,Z$. Every element $T$ of this set is a bounded linear operator from the  tangent space with the Banach norm to the tangent space with the Riemannian $g$-norm.
By the Banach-Steinhaus theorem, there exists a uniform bound for all such linear operators, i.e. $||R(X,Y)Z||_g\leq C$ for some positive constant $C>0$.  
\end{proof}
\begin{rem}
Proposition~\ref{wholecuraturetensor} implies that for any $N>0$ if there exists 
unit vectors $X_N, Y_N, Z_N\in T_x M$ such that $||R(X_N,Y_N)Z_N||_g>N$, then for any fixed pair $Y,Z$, there exists a
unit vector $X_N\in T_x M$ such that $||R(X_N,Y)Z||_g>N$. This can be done for fixed $X,Z$ or $X,Y$ as well.
\end{rem}
Now we can use a similar argument to show that the boundedness of  curvature at each point of a small enough neighborhood $U_x$ of $x$ is equivalent to the boundedness of curvature at its neighborhood for a Banach manifold.
\begin{prop}
Let $M$ be a Banach manifold, with  bounded curvatures for all $p$ in a neighborhood of $x\in M$. Then there exists a neighborhood $U_x$ such that the curvature is uniformly bounded in  $U_x$.     
\end{prop}
\begin{proof}
Firstly we take a neighborhood of $x$ homeomorphic to an open set of a Banach space where there is a trivialising map of the tangent bundle. Without loss of generality, it suffices to prove the assertion for $U_x$ being such an open set of a Banach space. Accordingly, we can define the curvature map under trivialization $T:U_x\times TU_x\to U_x\times TU_x$ via $T(X):=R(X,Y)Z$, which gives rise to a collection of bounded linear operators $T_p|_{Y_p, Z_p}:p\times TU_x\to p\times TU_x$ via $T_p|_{Y_p, Z_p}(X_p):=R_p(X_p,Y_p)Z_p$ for any fixed point $p\in U_x$ and any fixed corresponding vectors $Y_p,Z_p$. By the fact that the tangent space $TU_x$ is a Banach space, it follows by the Banach-Steinhaus Theorem that there exists a uniform bound for all such linear operators, i.e. $||R(X,Y)Z||_g\leq C$ for some positive constant $C>0$.
\end{proof}
\begin{rem}
 It remains unknown whether $R(X,Y)Z$ is bounded at every point of one of its neighborhood $U_x$ or even uniformly bounded on $U_x$ when we only know that $R(X,Y)Z$ is bounded at a point $x$. For a simple example, where this case is in fact true, one can consider a group of diffeomorphisms endowed with a right (left) invariant metric, in which case, the geometry on any other points is equivalent to the geometry on the identity.  
\end{rem}
In particular, if the Riemannian metric is strong on $M$, then its curvature tensors are bounded for any point $x\in M$, as shown in the following proposition. 
\begin{prop}[Theorem 3.3 in~\cite{bauer2021}]
Let $(M,g)$ be a strong Riemannian metric on $M$, then the curvature is bounded  at any point $x\in M$. In particular, since every finite-dimensional Riemannian metric is strong, the curvature is locally bounded for finite dimensional manifolds.
\end{prop}
\begin{proof}
 Since $g$ is strong, $M$ is a Hilbert manifold with respect to the Riemannian metric $g$, and thus every smooth function on $M$ is continuous. The curvature is a smooth section of the bundle of tensors over $M$. By continuity it is locally bounded on $M$, that is to say, there exists a constant $C>0$ and a neighborhood $U$ of $x$ such that $||R(X,Y)Z||_g\leq C$ whenever $X, Y, Z\in TU_x$ are normalized vectors.
\end{proof}
In general, the curvatures of a weak Riemannian metric may be locally unbounded at some point $x\in M$. In the following we present a simple example.
\begin{expl}[Example 13.4 in~\cite{bauer2021}]
For a sequence $m:\mathbb N\to (0,\infty)$ with $\operatorname{lim}_{i\to\infty}m_i=0$, we consider the weak Riemannian metric on $\mathfrak{l}^2$ given by
\begin{align*}
G_x(u,v) :=e^{-\frac{1}{2}||x||_{\mathfrak{l}^2}^2} ||u||_{\mathfrak{l}_m^2}:=e^{-\frac{1}{2}||x||_{\mathfrak{l}^2}^2}\sum_i m_i u_iv_i,\ \  x\in \mathfrak{l}^2, \ u,v\in T_x\mathfrak{l}^2=\mathfrak{l}^2.
\end{align*} 
The Christoffel symbol of this metric exists and is given by
\begin{align*}
\Gamma_x(u,v)=\frac{1}{2}\langle x,u\rangle_{\mathfrak{l}^2}v+\frac{1}{2}\langle x,v\rangle_{\mathfrak{l}^2}u-\frac{1}{2}\langle u,v\rangle_{\mathfrak{l}_m^2} m^{-1}x.   
\end{align*}
Consequently, the derivative of the Christoffel symbol satisfies
\begin{align*}
d\Gamma_x(v,w)=\frac{1}{2}\langle u,v\rangle_{\mathfrak{l}^2}w+\frac{1}{2}\langle u,w\rangle_{\mathfrak{l}^2}v-\frac{1}{2}\langle v,w\rangle_{\mathfrak{l}_m^2} m^{-1}u.   
\end{align*}
Then the curvature at $x=0$ satisfies
\begin{align*}
\langle R_0(u,v)v,u\rangle_{g_0}&=-\langle d\Gamma_0(u)(v,v),u\rangle_{g_0}+\langle d\Gamma_0(v)(u,v),u\rangle_{g_0}\\
&=-\langle u,v\rangle_{\mathfrak{l}^2}\langle v,u\rangle_{\mathfrak{l}_m^2}+\frac{1}{2}\langle u,u\rangle_{\mathfrak{l}^2}\langle v,v\rangle_{\mathfrak{l}_m^2}+\frac{1}{2}\langle v,v\rangle_{\mathfrak{l}^2}\langle u,u\rangle_{\mathfrak{l}_m^2}.   
\end{align*}
For $i\neq j$, the two normalized vectors $u:=m_i^{-1/2}e_i$ and $v:=m_j^{-1/2}e_j$ of the canonical vectors $e_i$ and $e_j$ are orthogonal with respect to $g_0=\langle\cdot,\cdot\rangle_{\mathfrak{l}_m^2}$ and satisfy
\begin{align*}
\langle R_0(u,v)v,u\rangle_{g_0}=m_i^{-1}+m_j^{-1}    
\end{align*}
Thus the sectional curvature is unbounded at $x=0$.
\end{expl}

\section{The existence of a local geodesic ball}\label{sec:geodesicball}
In later sections, we will need the existence of a geodesic ball around an arbitrary point to consider its Jacobi fields and conjugate points in the ball. Therefore, in this preparatory section, we first consider for which conditions we indeed have a local geodesic ball.

Let $M$ be a Banach manifold, and let $g$ be a Riemannian metric on $M$. It is clear that if $g$ is a strong metric, then a local geodesic ball must exist for any point in $M$. Unfortunately, this is not true for an arbitrarily chosen weak metric. The only thing we can obtain here is  the existence and uniqueness properties of local flows on a Banach manifold. This property and its proof can be found in~\cite{lang1962}, and we present it in the following lemma.
\begin{lem}\label{localflow}
Let $X$ be a smooth vector field on a manifold
$M$. If $M$ is a Banach manifold, then $X$ has a local flow. If $X$ has a local flow, then every
integral curve $c$ of $X$ satisfies $c(t) = \operatorname{Fl}^X_t (c(0))$. Thus, in this case there exists a unique
maximal flow.
\end{lem}
\begin{proof}
 The existence of a local flow on a Banach manifold can be verified in charts by using the Picard-Lindel\"{o}f theorem. Uniqueness follows from the following computation:
 \begin{align*}
\frac{d}{dt}\operatorname{Fl}^X (-t,c(t))&=-\frac{d}{ds}|_{s=-t}\operatorname{Fl}^X (s, c(t))+\frac{d}{ds}|_{s=t}\operatorname{Fl}^X (-t, c(s))\\
&=-\frac{d}{ds}|_{s=0}\operatorname{Fl}^X_{-t}\operatorname{Fl}^X (s, c(t))+T(\operatorname{Fl}^X_{-t})\dot c(t)\\
&=-T(\operatorname{Fl}^X_{-t}) X( c(t))+T(\operatorname{Fl}^X_{-t})X(c(t))=0,
 \end{align*}
 which implies that $\operatorname{Fl}^X (-t,c(t))=c(0)$ is a constant and therefore $c(t)=\operatorname{Fl}^X_t (c(0))$. 
\end{proof}
It is easily seen that for an arbitrary vector field, the infimum of the maximal existence of flow in all directions may vanish, which leads to the failure of a local geodesic ball.
In order to retain the existence of such a geodesic ball in infinite-dimensional geometry, we have to restrict ourselves to some specific cases. 

One way to overcome this difficulty is to expand the given Banach manifold $M$ to a larger one $\tilde M$. Then in some cases, the smooth geodesic neighborhood can be replaced by a weak geodesic ball.

\begin{expl}
Using the geometric interpretation of the Hunter-Saxton equation
as the Euler-Arnold equation for the geodesic flow on $\operatorname{Diff}(S^1)$ endowed with $\dot H^1$-metric, we are able to show that 
the maximal existence of time for such a geodesic flow is given by
\begin{equation}
T^*(u_0)=\frac{\pi}{2}+\operatorname{arctan}(\frac{1}{2}\operatorname{min}_{x\in S^1}u_{0x}(x))< \frac{\pi}{2},   
\end{equation}
where $u_0$ and $u_{0x}$ are the initial vector and its derivative. Therefore, the infimum of all positive numbers $T^*$ is zero, since the minimum of $u_{0x}$ can be as low as possible. However, we can extend the solution beyond its breaking time. In his paper~\cite{lenells2007}, Lenells considered the set $M^{\operatorname{AC}}$  of non-decreasing absolutely continuous functions $\varphi: [0,1]\to [0,1]$, and showed that the the geodesic flow in  $M^{\operatorname{AC}}$ is global as a weak solution to Euler-Arnold equation. 
\end{expl}
The other way is more subtle. We can endow the Riemannian  Banach manifold $M$ with extra structures, e.g., a Banach Lie group. Within this scope, we are able to show that the bounded covariant derivatives of curvatures lead to a local geodesic ball around the identity.  
The first step to accomplish this is to show that the Christoffel symbols with respect to $g$ can be chosen to be bounded within a ball by the boundedness of covariant derivatives.
\begin{lem}
 Let $M$ be a Banach Lie group endowed with a right invariant weak Riemannian metric with uniformly bounded covariant derivatives $\nabla^m R$ for all $m\geq0$. Then the Christoffel symbols with respect to a normal coordinate around $p$ are uniformly bounded along the geodesics.   
\end{lem}

\begin{proof}
Now we consider any geodesic $\gamma$ from the identity of the weak metric $g$, for which local existences are given by Lemma~\ref{localflow}. Since it is a normal coordinate at $p$, it is clear that the Taylor expansion of Christoffel symbols around $p$ is given by 
\begin{align*}
 \Gamma_{ij}^k(x)=-\frac{1}{3}\left(R_{kijl}(p)+R_{kijl}(p)\right)x_l+O(|x|^2),\ \  \forall i,j,k,  
\end{align*}
 where we can assume $\gamma'$ to be labeled by $l$ in the basis of the normal coordinate at $p$. Since the metric is right invariant, the curvature tensor $R$
has constant components in a right invariant frame. Therefore, all covariant derivatives $\nabla^m R$ also have constant components in this frame. Thus for $x$ in a small ball $B_\delta(p)$, there exists a constant $K>0$ such that 
\begin{align*}
 |\partial_l\Gamma_{ij}^k(x)|\leq \frac{2}{3}C+K|x|.  
\end{align*}
By a straightforward calculation, the Christoffel symbols in $B_\delta(p)$ are bounded by
\begin{align*}
|\Gamma_{ij}^k(\gamma(t))|\leq \frac{2}{3}C\delta+\frac{1}{2}K\delta^2.     
\end{align*}
\end{proof}

Next we use the Picard-Lindel\"{o}f theorem to give a uniform lower bound of the existence for any parametrized geodesic from the identity. 
\begin{thm}\label{existencegeodesicball}
 Let $M$ be any Banach Lie group  endowed with a right invariant weak Riemannian metric with bounded covariant derivatives $\nabla^m R$ for $m\geq0$. Then there exists a local geodesic ball around the identity element. 
\end{thm}
\begin{proof}
 It is easy to see that the prescribed geodesic equation can be written as 
 \begin{align*}
 \gamma'(t)=f(t,\gamma(t)), \ \ \gamma(0)=\operatorname{Id},   
 \end{align*}
 where $f$ is a matrix function of coordinate functions and its first derivatives with coefficients given by Christoffel symbols. Since Christoffel symbols are uniformly bounded within a small interval $[0,\delta]$ by the preceding lemma, $f$ is Lipschitz continuous in $\gamma$. Then an application of the Picard-Lindel\"{o}f Theorem shows that there exists a unique solution to the geodesic equation in $[0,\epsilon]$ for $\epsilon$ being a positive constant dependent only on $C$, 
 $K$ and $\delta$.
\end{proof}

\section{Estimate of the norm of the Jacobi field and exponential map}
In this section, we estimate the length of the Jacobi field around an arbitrary point and show that the length $||J||$ is positive in a small ball under the bounded curvature condition. Using this result, we are able to show that the exponential map is a local diffeomorphism in a small ball.

Firstly we prove a variation of the Sturm comparison theorem, and it can be seen as a version of the Rauch comparison theorem. For more details, one can refer to Chapter 10 in~\cite{docarmo1992}.

\begin{lem}\label{onedimthm}
A nonnegative function $f$ satisfies $f(0)=0, f'(0)=1$ and $f''(t)+Af(t)\geq0$ for some positive constant $A>0$. Then the first zero of $f$ does not occur before the first zero of $\tilde f(t)=\frac{1}{\sqrt{A}}\operatorname{sin}(\sqrt{A}t)$.   
\end{lem}
\begin{proof}
First note that $\tilde f(t)=\frac{1}{\sqrt{A}}\operatorname{sin}(\sqrt{A}t)$ is the unique solution to the differential equation $\tilde f''(t)+A\tilde f(t)=0$ under the given initial condition. It suffices to show that if $\tilde f(t)>0$ on $(0,t_0)$ and  $\tilde f(t_0)=0$, then $f(t)>0$ on $(0,t_0)$. By use of integral by parts, we obtain the following inequality
\begin{align*}
0<&\int_0^t\{\tilde f(f''+Af)-f(\tilde f''+A\tilde f) \}  dt
\\&=(\tilde ff'-f\tilde f')|_0^t.
\end{align*}
Note that the strict inequality on the left hand side holds since otherwise we have $f''+Af\equiv0$.
From the initial condition, $f(t)$ is positive in a neighborhood of zero. If $t_1$ is the smallest number in $(0,t_0)$ such that $f(t_1)=0$, then $f'(t_1)\leq0$, and we already have $\tilde f(t_1)>0$, which leads to a contradiction.
\end{proof}

Now we assume that the local geodesic ball exists and consider the Jacobi field $J$ with respect to the Riemannian metric $g$.
\begin{lem}\label{estimateofj}
Assume that all geodesics emanating from $p\in M$ exist for all initial vectors $v$ with $||v||\leq\delta$ for some $\delta>0$ and the curvature tensor $R$ is bounded.
Then there exists a small ball such that the length $||J||>0$ within this area.
\end{lem}
\begin{proof}
First by using the equation of the Jacobi field, we get an estimate on the second derivative of $J$
\begin{align*}
    ||\frac{D^2 J}{dt^2}||= ||R(\frac{d\gamma}{dt},J)\frac{d\gamma}{dt}||\leq A||J||,
\end{align*}
where the positive constant $A<\infty$ is the bound with respect to the curvature operator $R$. Next we consider the second derivative of $||J||$
\begin{align*}
\frac{d^2}{dt^2}||J||&=\frac{d}{dt}\frac{\langle\frac{DJ}{dt},J\rangle}{||J||}\\&=\frac{\langle\frac{DJ}{dt},\frac{DJ}{dt}\rangle||J||^2-\langle\frac{DJ}{dt},J\rangle^2}{||J||^3}+\frac{\langle\frac{D^2J}{dt^2},J\rangle}{||J||}\\&\geq -A||J||.
\end{align*}
The last inequality holds since the first term is positive by Cauchy-Schwarz. By above Lemma~\ref{onedimthm}, we obtain a lower bound of $||J||$, the function $\operatorname{sin}(\sqrt{A}t)$, which gets its nearest zero at $\frac{\pi}{\sqrt{A}}$. This gives a small ball for $||J||$ being uniformly positive.
\end{proof}

Before proving the main theorem of this section, we present a well-known result of the exponential map, and the proof can be found in~\cite{lang1962}.

\begin{lem}\label{expadjoint}
Let $\gamma:[0,1]\to M$ be the geodesic connecting $p$ and $q$ with initial vector $v_0$ and final vector $v_1$. If $w_0\in T_p M$ and $w_1\in T_qM$, then 
\begin{align*}
\langle (d\operatorname{exp}_p)_{v_0}(w_0), w_1\rangle=\langle w_0,(d\operatorname{exp}_q)_{-v_1}(w_1)\rangle.
\end{align*}
In particular, if $(d\operatorname{exp}_p)_{-v_1}$ is an injective, then 
$(d\operatorname{exp}_p)_{v_0}$ is also a surjective, and vice versa.
\end{lem}
\begin{proof}
We only prove the second assertion. For contradiction, we assume that $(d\operatorname{exp}_p)_{v_0}$ is not a surjective.
Then there is a nonzero element $w_1$ orthonormal to the whole image of $(d\operatorname{exp}_p)_{v_0}$. By the arbitrariness of vector $w_0$, $(d\operatorname{exp}_q)_{-v_1}(w_1)$ vanishes, which leads to $w_1=0$ since $(d\operatorname{exp}_q)_{-v_1}$ is an injective, a contradiction.
\end{proof}

Putting together all of the lemmas above, we are able to find a small normal ball $U$ of the zero section such that the differential $(d\operatorname{exp}_p)_v$ of the exponential map is a bijection at every point of $U$. 
\begin{prop}\label{normalball}
Let $M$ be a Banach manifold with locally bounded curvature and the lower-bound of maximal intervals of existence of geodesics is positive for any initial vectors with norm no more than $1$.  For some positive constant $\epsilon$ only pertaining to the upper bound of curvatures, there exists a normal ball $B(p,\epsilon)\subset T_p M$ of any point $p\in M$ such that the differential of the exponential map $(d\operatorname{exp}_p)_v$ is a bijection at every point $v\in B(p,\epsilon)$. This implies that the exponential map is a local diffeomorphism onto its image.   
\end{prop}
\begin{proof}
Let $\delta>0$ be the uniform lower-bound of the maximal interval of existence of geodesics, and some constant $C>0$ be the upper-bound of  curvatures. We observe that for some $0<\epsilon<\operatorname{min}(\frac{\pi}{\sqrt{C}},\delta)$, the length of the Jacobi field $||J(t)||>0$ for $t\in[0,\epsilon]$. Note that $\epsilon$ only depends on $C$, but not on the choice of the initial vector of the geodesic $v$ or the geodesic spray $w$. Now we consider the exponential map $\operatorname{exp}_p$ on a normal ball $B(p,\epsilon)$. Then the differential of the exponential map $(d\operatorname{exp}_p)_v$ in $B(p,\epsilon)$ is injective since every nonzero vector $w\in T_0(T_p M)$ will be transported to a nonzero vector in $T_v(T_p M)$. The surjectivity of $d\operatorname{exp}_p$ directly comes from the injectivity and Lemma~\ref{expadjoint}.  

Next we show that the exponential map is a local diffeomorphism on $B(p,\epsilon)$. Note that the differential $(d\operatorname{exp}_p)_v$ is a local bounded linear isomorphism from the above discussion. Then it is a local diffeomorphism by using the inverse function theorem for Banach manifolds.

Finally let $F$ denote the exponential map, then for each point $y\in\operatorname{exp}_p(B(p,\epsilon))$, there exists a small neighborhood $U_y$ of $y$ and a continuously differentiable map $G:U_y\to B(p,\epsilon)$ such that $F(G(z))=z$ for all $z\in U_y$. As $dF$ is an isomorphism at each point $v\in B(p,\epsilon)$, we can find a local inverse $G$ for any point $y\in\operatorname{exp}_p(B(p,\epsilon))$. Since this $G$ is the unique solution to the equation $F(a)=b$, $G$ is well defined globally on $\operatorname{exp}_p(B(p,\epsilon))$, which concludes the proof. 
\end{proof}

\section{The vanishing geodesic distance phenomenon}\label{sec:vanishing}

In this final section, we show that the locally unbounded curvature is a necessary requirement for vanishing geodesic distance phenomenon to occur in some specific cases. In other words, we need to show that the geodesic distance would separate points for locally bounded curvatures. We recall first the Gauss's lemma, note that it also holds for infinite dimensional setting.
\begin{lem}
Let $\gamma(t,s)=\operatorname{exp}_p(tv(s))$ be a family of geodesics from a fixed point $p$. Then we have $\langle\frac{\partial \gamma}{\partial t}, \frac{\partial \gamma}{\partial s}\rangle|_{(1,0)}=0$    
\end{lem}

Next we generalize Proposition 3.3.6 in Do Carmo's book~\cite{docarmo1992} and obtain the main theorem of the present article.
\begin{defi}[Property~B]
We say that a Banach manifold $M$ with Riemannian metric $g$ satisfies Property~B if it satisfies one of the following conditions:
\begin{itemize}
    \item [(1).] $M$ has a local geodesic ball around every point $p$ in $M$.
    \item [(2).] $(M,g)$ can be embedded in a larger Banach manifold 
    $(\tilde M,\tilde g)$
    which satisfies condition (1).
\end{itemize}

\end{defi}
\begin{thm}[Main Theorem]
Let $(M,g)$ be a Riemannian Banach manifold satisfying Property~B and having locally bounded curvatures. Then the geodesic distance of any two distinct points in $M$ is non-vanishing.   
\end{thm}
\begin{proof}
We introduce Property~B to ensure the existence of a local geodesic ball. Note that condition (2) is also sufficient since if the geodesic distance between two points $p$ and $q$ is non-vanishing in $\tilde M$, then the corresponding geodesic distance in the submanifold $M$ is larger, which is also non-vanishing.

So it remains to show that the geodesic distance indeed separates points when the existence of a local geodesic ball is guaranteed. Next we give a proof following the same steps in~\cite{docarmo1992}.

For any two distinct points $p$ and $q$, we consider the set of all smooth curves $c:[0,1]\to\mathbb R$ connecting $c(0)=p$ and $c(1)=q$.
Firstly we assume that $c([0,1])\subset \operatorname{exp}_p(B(p,\epsilon))$, where $B(p,\epsilon)$ is a normal ball in the tangent space, whose existence is guaranteed by Proposition~\ref{normalball}. Since $\operatorname{exp}_p$ is a local diffeomorphism, the curve $c(t)$, for $t\neq0$, can be written uniquely as $\operatorname{exp}_p(r(t)v(t))=\gamma(r(t),t)$ where $v(t)$ satisfies $||v(t)||=1$ and $r:(0,1]\to \mathbb R$ is a positive piecewise differentiable function. It follows that 
\begin{align*}
    \frac{dc}{dt}=\frac{\partial\gamma}{\partial t}r'(t)+\frac{\partial\gamma}{\partial s}
\end{align*}
By Gauss's lemma,
$\langle\frac{\partial\gamma}{\partial t}, \frac{\partial\gamma}{\partial s}\rangle=0$. Since $||\frac{\partial\gamma}{\partial t}||=1$, we have
\begin{align*}
    ||\frac{dc}{dt}||^2=|r'(t)|^2+||\frac{\partial\gamma}{\partial s}||^2\geq |r'(t)|^2
\end{align*}
Thus $\int_\delta^1 ||\frac{dc}{dt}||dt\geq\int_\delta^1 r'(t)dt=r(1)-r(\delta)$.

Taking $\delta\to0$, the geodesic distance between $p$ and $q$ has a lower bound $r(1)>0$.
If $c([0,1])$ is not contained in $B(p,\epsilon)$, consider the first point $t_0$ for which $c(t_0)\in \partial B(p,\epsilon)$. Then the length of $c$, $l(c)\geq l_{[0,t_0]}(c)\geq\epsilon$, which implies the geodesic distance  between $p$ and $q$ has a lower bound $\epsilon>0$.
\end{proof}

\appendix
\section{Semi-boundedness of curvature tensors}\label{sec:appendix}
In Proposition~\ref{wholecuraturetensor}, we see that the curvature tensor  is bounded if it is bounded when two of the vectors are fixed. Now if we loosen the restriction, that is to say, the two fixed vectors are identical, then we obtain a weaker boundedness condition on curvature tensors, which is denoted by semi-boundedness here.

\begin{defi}[Semi-boundedness of curvature]
A  Riemmannian manifold $(M,g)$ is said to have a semi-bounded curvature at $x\in M$ if for any fixed vector field $X\in TM$,
there exists a constant $C_X>0$ such that the curvature tensor $||R(X,Y)X||_{g_x}\leq C_X$ for any unit vectors $Y_x\in T_x M$. Note that the location of $Y$ in $R(X,Y)X$ does not change the semi-boundedness, since $R(X,Y)X=-R(Y,X)X$.  
\end{defi}
It is clear to see that the curvature is bounded then it is semi-bounded. For the converse, we need the uniform boundness of $C_Y$.
\begin{prop}
If there exists a constant $C>0$ such that $C_X\leq C$ for all unit vectors $X$, then the semi-boundedness implies the boundness on curvatures.    
\end{prop}
\begin{proof}
By the Cauchy-Schwarz inequality, we have for any unit vectors $X$ and $Y$
\begin{align*}
||\langle R(X,Y)X,Y\rangle||\leq ||R(X,Y)X||||Y|| \leq C_X\leq C.   
\end{align*}
Then the result follows by Proposition~\ref{seccurv}.
\end{proof}

This proposition implies that if a Riemannian metric on a infinite dimensional manifold $M$ has semi-bounded curvatures but not bounded curvatures, then there exists a sequence of vectors 
$X$ such that $C_X$ goes to $\infty$. It is also important to note that there indeed exists such curvature tensors. For instance, one can consider the right-invariant weak $L^2$-metric on the volumorphism group $\operatorname{Diff}_\mu^s(M)$ of Sobolev class $H^s$ of $M$ for $s>n/2+1$; see~\cite{misiolek2010}.
\begin{lem}[Lemma 8.1 in~\cite{misiolek2010}]
For any vectors $X,Y\in T_e\operatorname{Diff}_\mu^s(M)$, we have the following estimate
\begin{align*}
||R(X,Y)X||_{H^\sigma}\leq C||X||^2_{H^s}||Y||_{H^\sigma} \ \ (0\leq\sigma\leq s),   
\end{align*}
where the constant $C$ depends on  $s>n/2+1$.    
\end{lem}
It is easy to see that the semi-bound of curvature tensors is $C_X=C||X||^2_{H^s}/||X||^2_{H^\sigma}$, and it can go to infinity if the vectors $X$ are chosen appropriately for $\sigma<s$. 

Although for semi-bounded curvatures we fall short of good properties regarding its exponential map and  local geodesic ball, 
it is still possible to say something about its Jacobi field.
\begin{prop}
Let $M$ be a Riemannian manifold with semi-bounded curvatures. Then conjugate points on any finite geodesic segment $\gamma$ are isolated, and the geodesic contains at most finitely many points which are conjugate to $p=\gamma(0)$ along $\gamma$.
\end{prop}
\begin{proof}
This can be easily seen from the proof of Lemma~\ref{estimateofj}. 
The distance between two conjugate points is larger than $\frac{\pi}{\sqrt{N_X}}$ by an application of the comparison theorem.
\end{proof}

\subsection*{Acknowledgements}
This article is motivated by amounts of communication with the author's advisor M. Bauer, and the author is grateful to his advisor for his informative suggestions and patience.




\bibliographystyle{abbrv}

\end{document}